\magnification\magstephalf
\def\meno{\medskip\noindent}
\def\proof{\noindent{\bf Proof.}\quad}
\def\pfbox
  {\hbox{\hskip 3pt\lower0pt\vbox{\hrule
  \hbox to 6pt{\vrule height 6pt\hfill\vrule}
  \hrule}}\hskip3pt}
\def\bib[#1] {\par\noindent\hangindent 20pt\hbox to20pt{[#1]\hfil}}
\def\AW{Addison\kern.1em--Wesley}

\def\ra{\rightarrow}
\def\tta{\tt a}
\def\ttb{\tt b}
\def\ttc{\tt c}
\def\ttd{\tt d}
\def\tte{\tt e}
\def\ttf{\tt f}
\def\ttg{\tt g}
\def\tth{\tt h}
\def\tti{\tt i}
\def\ttj{\tt j}
\def\ttk{\tt k}
\def\ttl{\tt l}
\def\ttm{\tt m}
\def\ttn{\tt n}
\def\tto{\tt o}
\def\ttp{\tt p}
\def\ttq{\tt q}
\def\ttr{\tt r}
\def\tts{\tt s}

\def\ttu{\tt u}
\def\ttv{\tt v}
\def\ttw{\tt w}
\def\ttx{\tt x}
\def\tty{\tt y}
\def\ttz{\tt z}
\def\ttA{\tt A}
\def\ttB{\tt B}
\def\ttC{\tt C}
\def\ttD{\tt D}
\def\ttE{\tt E}
\def\ttF{\tt F}
\def\ttG{\tt G}
\def\ttH{\tt H}
\def\ttI{\tt I}
\def\ttJ{\tt J}
\def\ttK{\tt K}
\def\ttL{\tt L}
\def\ttM{\tt M}
\def\ttN{\tt N}
\def\ttO{\tt O}
\def\ttP{\tt P}
\def\ttQ{\tt Q}
\def\ttR{\tt R}
\def\ttS{\tt S}
\def\ttT{\tt T}
\def\ttU{\tt U}
\def\ttV{\tt V}
\def\ttW{\tt W}
\def\ttX{\tt X}
\def\ttY{\tt Y}
\def\ttZ{\tt Z}

\def\Ga{1}
\def\GR{2}
\def\GK{3}
\def\ii{4}
\def\Ry{5}

\centerline{\bf The Knowlton-Graham Partition Problem}
\smallskip
\centerline{Donald E. Knuth}
\centerline{\sl Computer Science Department}
\centerline{\sl Stanford University}
\centerline{\sl Stanford, CA 94305--2140}

\bigskip
{\narrower\smallskip\noindent
{\bf Abstract.}\def\ttz{\tt z}\def\ttz{\tt z}
A set partition technique that is useful for identifying wires in cables can be
recast in the language of 0--1 matrices, thereby resolving an open problem
stated by R.~L. Graham in Volume~1 of this journal. The proof involves a
construction of 0--1 matrices having row and column sums without gaps.
\smallskip}

{\baselineskip14pt

\meno
A  long cable contains $n$ indistinguishable wires. Two people, one at each
end, want to label the wires consistently so that both ends of each wire
receive the same label. An interesting way to achieve this was proposed by
K.~C. Knowlton
[\GK]:
Partition $\{1,\ldots,n\}$ into disjoint sets in two ways $A_1,\ldots,A_p$ and
$B_1,\ldots,B_q$, subject to the condition that at most one element appears
both in an $A$~set of cardinality~$j$ and in a $B$~set of cardinality~$k$, for
each $j$ and~$k$. We can then use the coordinates $(j,k)$ to identify each
element. R.~L. Graham
[\GR]
proved that such partitioning schemes exist if and only if $n\neq 2$, 5, or~9.

By restating the problem in terms of 0--1 matrices, it is possible to prove
Graham's theorem more simply, and to sharpen the results of
[\GR].

\proclaim
Lemma 1. Knowlton-Graham partitions for $n$ exist if and only if there is a
matrix of 0s and~1s having row sums $(r_1,\ldots,r_m)$ and column sums
$(c_1,\ldots,c_m)$ such that $r_j$ and~$c_j$ are multiples of~$j$ and
$r_1+\cdots +r_m=c_1+\cdots +c_m=n$.

\proof
If $A_1,\ldots,A_p$ and $B_1,\ldots,B_q$ are partitions of $\{1,\ldots,n\}$
with the Knowlton-Graham property, let $a_{jk}$ be the number of elements that
appear in an $A$~set of cardinality~$j$ and a $B$~set of cardinality~$k$. Then
$a_{jk}$ is 0 or~1; and $r_j=\sum_ka_{jk}$ is $j$~times the number of $A$~sets
of cardinality~$j$, while $c_k=\sum_ja_{jk}$ is $k$~times the number of
$B$~sets of cardinality~$k$.

Conversely, given such a matrix, we can use its rows
to define $A_1,\ldots,A_p$ such that each~1 in row~$j$ is in an $A$~set
of cardinality~$j$; similarly its columns define $B_1,\ldots,B_q$
such that each~1 in column~$k$ is in a $B$~set of cardinality~$k$. \ \pfbox

For example, the symmetric matrix
$$
\advance\baselineskip-2pt
\vcenter{\halign{\hfil#\hfil\ &\hfil#\hfil\ &\hfil#\hfil\ %
&\hfil#\hfil\ &\hfil#\hfil\ &\hfil#\hfil\cr
0&1&0&0&0&1\cr
1&1&1&1&1&1\cr
0&1&0&0&1&1\cr
0&1&0&1&1&1\cr
0&1&1&1&1&1\cr
1&1&1&1&1&1\cr}}$$
has row and column sums $(2,6,3,4,5,6)$ that satisfy the divisibility condition
and sum to~26. To identify 26~wires, we can associate the~1s with arbitrary
labels $\{\tta,\ldots,\ttz\}$,
$$
\advance\baselineskip-2pt
\vcenter{\halign{\hfil{\tt#}\hfil\ &\hfil{\tt#}\hfil\ %
&\hfil{\tt#}\hfil\ &\hfil{\tt#}\hfil\ &\hfil{\tt#}\hfil\ %
&\hfil{\tt#}\hfil\cr
.&a&.&.&.&b\cr
c&d&e&f&g&h\cr
.&i&.&.&j&k\cr
.&l&.&m&n&o\cr
.&p&q&r&s&t\cr
u&v&w&x&y&z\cr}}$$
The person at one end of the cable labels the wires with
$\{\tta,\ldots,\ttz\}$ arbitrarily and makes connections so that each element
of row~$j$ is connected to exactly $j-1$ other elements of its row; for
example, the connected components might be $A_1,\ldots,A_9=\{\tta\},\{\ttb\},
\{\ttc,\ttd\},\{\tte,\ttf\},\{\ttg,\tth\},\{\tti,\ttj,\ttk\},
\{\ttl,\ttm,\ttn,\tto\}$,
$\{\ttp,\ttq,\ttr,\tts\},\{\ttu,\ttv,\ttw,\ttx,\tty,
\ttz\}$. The person at the other end now uses properties of conductivity
to tell what row each wire belongs~to. 
The wires at that end can then be labeled
$\{\ttA,\ldots,\ttZ\}$ in such a way that
$\{\ttA,\ttB\}=\{\tta,\ttb\}$, $\{\ttC,\ttD,\ttE,\ttF,\ttG,\ttH\}
=\{\ttc,\ttd,\tte,\ttf,\ttg,\tth\}$, etc. Now the wires at the first end 
 are disconnected, while at the other end they are connected so that each
element of column~$k$ is connected to exactly $k-1$ other elements of its
column. For example, the connected components might now be $B_1,\ldots,B_9
=\{\ttC\},\{\ttU\},\{\ttA,\ttD\},\{\ttI,\ttL\},\{\ttP,\ttV\},
\{\ttE,\ttQ,\ttW\},\{\ttF,\ttM,\ttR,\ttX\},
\{\ttG,\ttJ,\ttN,\ttS,\ttY\}$,\break
$\{\ttB,\ttH,\ttK,\ttO,\ttT,\ttZ\}$. 
Once this has been done, the people at both ends of the cable
can give unique coordinates $(j,k)$ to each wire, knowing its row and column,

Knowlton-Graham partitions are said to have order $m$ if the largest
cardinality of $A_1,\ldots, A_p$, $B_1,\ldots,B_q$ is~$m$. 

\proclaim
Theorem 1. {\rm (Graham).}
Knowlton-Graham partitions of\/ $n$ having order\/ $m$ are possible only if\/
${m+1\choose 2}\leq n\leq J(m)$, where
$$J(m)=\sum_{j=1}^mj\lfloor m/j\rfloor\,.$$

\proof
By Lemma 1, Knowlton-Graham partitions of order $m$ imply the existence of an
$m\times m$ matrix of 0s and~1s having row sums $(r_1,\ldots,r_m)$ and column
sums $(c_1,\ldots,c_m)$, where both $r_j$ and~$c_j$ are multiples of~$j$ for
$1\leq j\leq m$, and where $r_m+c_m>0$. Clearly $r_j\leq m$; so $r_j$ is at
most $j\lfloor m/j\rfloor$, the largest multiple of~$j$ that does not
exceed~$m$. This establishes the upper bound $J(m)$.

If $r_m>0$, we must have $r_m=m$; this implies $c_j>0$ for all~$j$, hence
$c_m=m$. Similarly, $c_m>0$ implies that $r_m=c_m=m$. So we must have $r_j>0$
for all~$j$, hence $r_j\geq j$ for all~$j$, hence $n=\sum_{j=1}^m r_j\geq 
\sum_{j=1}^m j={m+1\choose 2}$. \ \pfbox

\medskip
When $m=1,2,3,4$, Theorem 1 says that $1\leq n\leq 1$, $3\leq n\leq 4$, $6\leq
n\leq 8$, and $10\leq n\leq 15$, respectively; this explains why the values
$n=2,5,9$ are impossible. For $m\geq 4$ we have $J(m)\geq {m+2\choose 2}$, so
there are no more gaps. In fact, as $m\ra\infty$ we have
$$J(m)=m^2-\sum_{j=1}^m(m\bmod j)={\pi^2\over 12}\,m^2+O(m\log m)$$
(see [\ii, Eq.\ 4.5.3--21]); 
therefore
$J(m)\left/{m+2\choose 2}\right.$ approaches the limiting value
$\pi^2\!/6\approx 1.64$.

The main purpose of this note is to prove the converse of Theorem~1, namely
that Knowlton-Graham partitions of order~$m$ do exist for all~$n$ in the 
range  ${m+1\choose 2}\leq n\leq J(m)$. This question was left open in 
[\GR],
where Graham observed that it was not sufficient simply to represent~$n$ in the
form $r_1+\cdots +r_m$ where each $r_j$ is a positive multiple of~$j$. For
example, there is no 0--1 matrix having row sums $(1,6,6,4,5,6)$. If
there were, we would necessarily have $c_1\ge4$, $c_2\ge4$, $c_3=6$,
$c_4=4$, $c_5=5$, $c_6=6$, and $c_1+\cdots+c_6>r_1+\cdots+r_6$.

Gale [\Ga] and Ryser [\Ry] independently found an elegant necessary and
sufficient condition for the existence of 0--1 matrices having given row and
column sums, but their theorem does not seem to lead easily to the result
needed here. Instead, we can use a direct recursive construction.

\proclaim
Lemma 2. Let $(r_1,\ldots,r_m)$ and $(c_1,\ldots,c_m)$ be integers satisfying
the conditions
$$\displaylines{\hfil r_1+\cdots +r_m=c_1+\cdots +c_m\,,\hfil\cr
\noalign{\smallskip}
\hfil m\geq r_1\geq \cdots \geq r_m\geq 0\,,\quad
m\geq c_1\geq \cdots\geq c_m\geq 0\,,\hfil\cr
\noalign{\smallskip}
\hfil r_{j+1}\geq r_j-1\quad{\rm and}\quad c_{j+1}\geq c_j-1\quad
{\rm for}\quad 1\leq j<m\,.\hfil\cr}$$
Then there exists an $m\times m$ matrix of 0s and 1s having row sums
$(r_1,\ldots,r_m)$ and column sums $(c_1,\ldots,c_m)$.

\proof
This is obvious when $m=1$, so we may assume inductively that $m>1$. Let
$p=r_1$ and $q=c_1$, and consider the numbers
$$\eqalign{(r'_1,\ldots,r'_{m-1})&=(r_2-1,\ldots,r_q-1,r_{q+1},\ldots,r_m)\cr
\noalign{\smallskip}
(c'_1,\ldots,c'_{m-1})&=(c_2-1,\ldots,c_p-1,c_{p+1},\ldots,c_m)\,.\cr}$$
The lemma will be proved if we construct an $(m-1)\times (m-1)$ matrix
of~0s and~1s having row sums $(r'_1,\ldots,r'_{m+1})$ and column sums
$(c'_1,\ldots,c'_{m-1})$, because we can achieve the desired result by
appending a new first row and a new first column. In fact it suffices, by row
and column permutations, to construct a 0--1 matrix with row and column sums
equal to the numbers $(r''_1,\ldots,r''_{m-1})$ and $(c''_1,\ldots,c''_{m-1})$
obtained by sorting $(r'_1,\ldots,r'_{m-1})$ and $(c'_1,\ldots,c'_{m-1})$ into
nonincreasing order.

Since $r''_1+\cdots +r''_{m-1}=r_1+\cdots +r_m-p-q+1=c_1+\cdots +c_m-p-q+1
=c''_1+\cdots+c''_{m-1}$, we can use the induction hypothesis if we verify
that $r''_1\leq m-1$, $r''_{m-1}\geq 0$, $r''_{j+1}\geq r''_j-1$; the similar
inequalities for $(c''_1,\ldots,c''_{m-1})$ follow by symmetry.

Suppose $r''_1=m$; this implies $r_{q+1}=m$ and $q<m$. Therefore
$r_q=\cdots=r_2=r_1=m$, and we have $(q+1)m\leq r_1+\cdots +r_m=c_1+\cdots
+c_m\leq m\,c_1=qm$, a~contradiction.

Suppose $r''_{m-1}<0$; this implies $r_q=0$. Therefore $(q-1)+\cdots +1+0\geq
r_1+\cdots +r_m=c_1+\cdots +c_m\geq q+(q-1)+\cdots +1$, another contradiction.

Suppose finally that $r''_{j+1}<r''_j-1$. This could happen only if
$r''_j=r_k$ and $r''_{j+1}=r_l-1$ for some~$k$ and~$l$ with $r_k>r_l$.
But we would not decrease $r_l$ unless we had also decreased $r_k$. \ \pfbox

\medskip
The  construction of Lemma 2 produces a symmetric matrix when
$(r_1,\ldots,r_m)=(c_1,\ldots,c_m)$. Let's say that Knowlton-Graham partitions
are {\it symmetric\/} if they correspond to a symmetric matrix. We are now
ready to prove the main result.

\proclaim
Theorem 2. Symmetric Knowlton-Graham partitions of $n$ having order~$m$ exist
whenever ${m+1\choose 2}\leq n\leq J(m)$.

\proof
When $n$ is in the stated range but not equal to ${m+1\choose 2}$, there is a
number $s\leq m/2$ such that we can write $n=t_1+\cdots +t_m$ where
$$\vcenter{\halign{\hfil $#\;$&$#$\hfil\quad&#\hfil\cr
t_j&=j\,,&for $s<j\leq m\,;$\cr
t_s&=ks&for some $k$, $1< k\leq \lfloor m/s\rfloor\,,$\cr
t_j&=j\lfloor m/j\rfloor\,,&for $1<j<s\,;$\cr
m-s&<t_1\leq m\,,&if $s>1$.\cr}}$$
When $n=J(m)$, this is true with $s=\lfloor m/2\rfloor$, $k=\lfloor
m/s\rfloor$, and $t_1=m$. Otherwise we can find such a representation by
first
representing $n+1$ and subtracting~1 from~$t_1$; then if $s>1$ and $t_1=m-s$,
we replace~$t_1$ by~$m$ and subtract~$s$ from~$t_s$; finally, if $t_s=s$, we
decrease~$s$ by~1.

The remaining case $n={m+1\choose 2}$ is simpler because we can write
$n=t_1+\cdots+t_m$ where $t_j=j$ for all~$j$. This is a representation of
essentially the same form but with $s=0$.

Notice that $t_j$ is a multiple of $j$, for $1\leq j\leq m$. We can also verify
that the set $\{t_1,\ldots,t_m\}$ consists simply of the consecutive elements
$\{t_{s+1},\ldots,t_m\}=\{s+1,\ldots,m\}$. For we have $t_s>s$; and
$t_j>m-s\geq s$ for all $j<s$, because $j\lfloor m/j\rfloor =m-(m\bmod j)$.

Let $(r_1,\ldots,r_m)$ and $(c_1,\ldots,c_m)$ be the numbers $(t_1,\ldots,t_m)$
sorted into nonincreasing order. Lemma~2 tells us how to construct a symmetric
0--1 matrix having these row and column sums. After an appropriate permutation
of rows and columns, the row and column sums can be made equal to
$(t_1,\ldots,t_m)$; and this yields Knowlton-Graham partitions, by Lemma~1.
\ \pfbox

}

\bigskip
\centerline{References}

\smallskip
\bib
[\Ga]
David Gale, ``A theorem on flows in networks,'' {\sl Pacific Journal of
Mathematics\/ \bf 7} (1957), 1073--1082.

\smallskip
\bib
[\GR]
R. L. Graham, ``On partitions of a finite set,'' {\sl Journal of Combinatorial
Theory\/ \bf 1} (1966), 215--223.

\smallskip
\bib
[\GK]
Ronald L. Graham and Kenneth C. Knowlton, ``Method of identifying conductors in
a cable by establishing conductor connection groupings at both ends of the
cable,'' U.S. Patent 3,369,177 (13~Feb 1968).

\smallskip
\bib
[\ii]
Donald E. Knuth, {\sl Seminumerical Algorithms}, second edition
(Reading, Massachusetts: \AW, 1981).

\smallskip
\bib
[\Ry]
H. J. Ryser, ``Combinatorial properties of matrices of zeros and ones,'' 
{\sl Canadian Journal of Mathematics\/ \bf 9} (1957), 371--377.

\bye